\documentclass[a4paper,12pt,twoside]{article}
\usepackage{latexsym, amsmath, amsfonts, amscd, amssymb, verbatim}
\usepackage{color}
\usepackage{lscape}

\def\tto{\;{\lower 1pt \hbox{$\rightarrow$}}\kern -10pt
	\hbox{\raise 2pt \hbox{$\rightarrow$}}\;}

\begin{document}
	\pagestyle{myheadings}
	
	\newtheorem{Theorem}{Theorem}[section]
	\newtheorem{Proposition}[Theorem]{Proposition}
	\newtheorem{Remark}[Theorem]{Remark}
	\newtheorem{Lemma}[Theorem]{Lemma}
	\newtheorem{Corollary}[Theorem]{Corollary}
	\newtheorem{Definition}[Theorem]{Definition}
	\newtheorem{Example}[Theorem]{Example}
	\renewcommand{\theequation}{\thesection.\arabic{equation}}
	\normalsize
	\setcounter{equation}{0}
	\title{\bf On Some Incremental Algorithms for the Minimum Sum-of-Squares Clustering Problem. Part 1: Ordin and Bagirov's Incremental Algorithm\footnote{This work was supported by National Foundation for Science $\&$ Technology Development (Vietnam) under grant number 101.01-2018.308 and China Medical University (Taichung, Taiwan). The second author was partially supported by the Grant MOST 105-2115-M-039-002-MY3. We would like to thank Dr. Nguyen Van Tuyen for useful discussions on an example.}}
	
	\medskip
	\author{Tran Hung Cuong\footnote{Department of Computer Science, Faculty of Information Technology, Hanoi University of  Industry, Minh Khai Ward, North Tu Liem District, Hanoi, Vietnam; email: tranhungcuong@haui.edu.vn, tranhungcuonghaui@gmail.com.},\ \, Jen-Chih Yao\footnote{Center for General Education, China Medical University, Taichung 40402, Taiwan; email: yaojc@mail.cmu.edu.tw.},\ \, Nguyen Dong Yen\footnote{Institute of Mathematics, Vietnam Academy of
			Science and Technology, 18 Hoang Quoc Viet, Hanoi 10307, Vietnam;
			email: ndyen@math.ac.vn.}}
	\maketitle
	\date{}

\medskip
\begin{quote}
\noindent {\bf Abstract.}  Solution methods for the minimum sum-of-squares clustering (MSSC) problem are analyzed and developed in this paper. Based on the DCA (\textbf{D}ifference-of-\textbf{C}onvex functions \textbf{A}lgorithms) in DC programming and recently established qualitative properties of the MSSC problem \cite{CYY_p1}, we suggest several improvements of the incremental algorithms of Ordin and Bagirov \cite{Ordin_Bagirov_2015} and of Bagirov \cite{Bagirov_2014}. Properties of the new algorithms are obtained and preliminary numerical tests of those on real-world databases are shown. Finite convergence, convergence, and the rate of convergence of solution methods for the MSSC problem are presented for the first time in our paper. This Part 1 is devoted to the incremental heuristic clustering algorithm of Ordin and Bagirov and the modified version proposed herein.


\end{quote}

\section{Introduction}
\setcounter{equation}{0}

Data clustering is applied widely in image segmentation, information retrieval, pattern recognition, pattern classification, network analysis, vector quantization and data compression, data mining and knowledge discovery business, document clustering and image processing (see, e.g., \cite[p.~32]{Aggarwal 2014} and \cite{Kumar 2017}). 

\medskip
Among different clustering criteria, \textit{the \textbf{M}inimum \textbf{S}um-of-\textbf{S}quares \textbf{C}lustering} (MSSC for short) criterion is one of the most used \cite{Bock_1998,Brusco 2006,Costa 2017,Du Merle 2000,Kumar 2017,PhamDinh_LeThi_2009,Peng 2005,Shereli 2005}. Accepting this criterion, one tries to make the sum of the squared Euclidean distances from each data point to the centroid of its cluster as small as possible. To solve the \textit{MSSC problem} means to divide a finite data set into a given number of disjoint clusters so that the just mentioned sum is minimal. 

\medskip
 It is well known that qualitative properties of an optimization problem are very helpful for its numerical solution. In addition to the fundamental necessary optimality condition given by Ordin and Bagirov \cite{Ordin_Bagirov_2015}, a series of qualitative properties of the MSSC problem have been obtained in our recent paper \cite{CYY_p1}. Namely, we have proved that the MSSC problem always has a global solution and, under a mild condition, the global solution set is finite and the components of each global solution can be computed by an explicit formula. Besides, complete characterizations of the nontrivial local solutions of the MSSC problem have been established. Moreover, we have analyzed the changes of the optimal value, the global solution set, and the local solution set of the  MSSC problem with respect to small changes in the data set. We have shown that the optimal value function is locally Lipschitz, the global solution map is locally upper Lipschitz, and the local solution map has the Aubin property, provided that the original data points are pairwise distinct. 

\medskip
 There are many algorithms to solve the MSSC problem (see, e.g., \cite{Bagirov_2008,Bagirov_2014,Bagirov_Rubinov_2003,BUW_2011,Bagirov_Yearwood_2006,Ordin_Bagirov_2015,Xie_2011}, and the references therein). Since it is a NP-hard problem \cite{Aloise 2009,Mahajan_NPHard_2009} when either the number of the data features or the number of the clusters is a part of the input, the fact that the existing algorithms can give at most some local solutions is understandable. 
 
 \medskip
\textit{The $k$-means clustering algorithm} (see, e.g., \cite{Aggarwal 2014}, \cite{Jain 2010}, \cite{Kantardzic 2011}, and \cite{MacQueen 1967}) is the best known solution method for the MSSC problem. To improve its effectiveness, the so-called \textit{global $k$-means clustering algorithms} have been proposed \cite{Bagirov_2008,BUW_2011,HNCM_2005,Lai_Huang_2010,LVV_2003,Xie_2011}.
 
  \medskip
 Since the quality of the computation results  greatly depends on the starting points, it is reasonable to look for good starting points.  The DCA (Difference-of Convex-functions Algorithms), which has been applied to the MSSC problem in \cite{Bagirov_2014,PhamDinh_LeThi_2009}, can be used for the purpose. 
 
 \medskip
 One calls a clustering algorithm \textit{incremental} if the number of the clusters increases step by step. As noted in \cite[p.~345]{Ordin_Bagirov_2015}, the available numerical results demonstrate that \textit{incremental clustering algorithms} (see, e.g., \cite{Bagirov_2008,Ordin_Bagirov_2015,HNCM_2005,Lai_Huang_2010}) are efficient for dealing with large data sets.

\medskip
In this paper, we are interested in analyzing and developing the \textit{incremental heuristic clustering algorithm} of Ordin and Bagirov~\cite{Ordin_Bagirov_2015} and the \textit{incremental DC clustering algorithm} of Bagirov \cite{Bagirov_2014}. By constructing some concrete MSSC problems with small data sets, we will show how these algorithms work. It turns out that, due to the exact stopping criterion, the computation by the second algorithm may not stop. We propose three modified versions for this algorithm. As concerning the first algorithm, a modified version with a more reasonable computation procedure is suggested.

\medskip
To the best our knowledge, finite convergence, convergence, and the rate of convergence of solution methods for the MSSC problem are given for the first time in the present paper.  

\medskip
The paper is split into two parts. This Part 1 is devoted to the incremental heuristic clustering algorithm of Ordin and Bagirov and the modified version proposed herein. Among other things, several fundamental properties of the latter are established. Part 2 presents the incremental DC clustering algorithm of Bagirov and the three modified versions we suggest for it. Two convergence theorems and two theorems about the $Q$-linear convergence rate of the first modified version of Bagirov's algorithm will be obtained by some delicate arguments. Preliminary numerical tests of the above-mentioned algorithms on real-world databases are shown.

\medskip
The organization of Part 1 is as follows. In Section~2, the MSSC clustering problem and its basic qualitative properties, which have been obtained in \cite{CYY_p1}, are recalled. The $k$-means algorithm is also presented in Section 2. The incremental heuristic clustering algorithm of Ordin and Bagirov~\cite{Ordin_Bagirov_2015} and a modified version are described and analyzed in Section~3, where two theorems on the behavior of the latter are proved.

\section{Basic Properties of the Clustering Problem}
\setcounter{equation}{0}

Let $A = \{{a}^1,...,{a}^m\}$ be a finite set of points (representing the data points to be grouped) in the $n$-dimensional Euclidean space $\mathbb R^n$ equipped with the scalar product $\langle x,y\rangle =\displaystyle\sum_{i=1}^n x_iy_i$ and the norm $\|x\|=\left(\displaystyle\sum_{i=1}^{n}x_i^2\right)^{1/2}$. The open ball (resp., closed ball) with center $a\in\mathbb R^n$ and radius $\varepsilon>0$ will be denoted by $B(a,\varepsilon)$ (resp., $\bar B(a,\varepsilon)$). For a subset $\Omega$ of an  Euclidean space, its \textit{convex hull} is denoted by ${\rm co}\,\Omega$. Given a positive integer $k$ with $k\leq m$, one wants to partition $A$ into disjoint   subsets $A^1,\dots, A^k,$ called \textit{clusters}, such that a \textit{clustering criterion} is optimized. 

\medskip
If one associates to each cluster $A^j$ a \textit{center} (or \textit{centroid}), denoted by ${x}^j\in \mathbb R^n$, then the following well-known \textit{variance} or \textit{SSQ (Sum-of-Squares) clustering criterion} (see, e.g., \cite[p.~266]{Bock_1998})
\begin{equation*}
	\psi(x,\alpha):=\frac{1}{m}\sum_{i=1}^{m}\left(\sum_{j=1} ^{k}\alpha_{ij}\|{a}^i-{x}^j\|^2\right)\ \;\longrightarrow\ \;\min,
\end{equation*}
where $\alpha_{ij}=1$  if ${a}^i\in A^j$ and $\alpha_{ij}=0$ otherwise, is used. Thus, the above partitioning problem can be formulated as the constrained optimization problem
\begin{equation}\label{basic_clustering_problem}\begin{array}{rl}
		\min\Big\{\psi(x,\alpha)\; \mid &\in\mathbb R^{n\times k},\ \mathbb{\alpha}=(\alpha_{ij})\in \mathbb{R}^{m\times k},\ \alpha_{ij}\in \{0, 1\},\\ & \displaystyle\sum_{j=1}^{k}\alpha_{ij}=1,\ i=1,\dots,m,\ j=1,\dots, k\Big\},
	\end{array}
\end{equation} 
where the centroid system $x=({x}^1,\dots,{x}^k)$ and the incident matrix $\mathbb{\alpha}=(\alpha_{ij})$ are to be found.

\medskip
Since \eqref{basic_clustering_problem} is a difficult \textit{mixed integer programming problem}, instead of it one usually considers (see, e.g., \cite[p.~344]{Ordin_Bagirov_2015}) next \textit{unconstrained nonsmooth nonconvex optimization problem}
\begin{equation}\label{DC_clustering_problem}
	\min\Big\{f(x):=\frac{1}{m}\sum_{i=1}^{m}\left(\min_{j=1,\dots,k} \|{a}^i-{x}^j\|^2\right)\, \mid\, x=({x}^1,\dots,{x}^k)\in\mathbb R^{n\times k}\Big\}.
\end{equation} 

Both models \eqref{basic_clustering_problem} and \eqref{DC_clustering_problem} are referred to as \textit{the minimum sum-of-squares clustering problem} (the MSSC problem). As the decision variables of  \eqref{basic_clustering_problem} and \eqref{DC_clustering_problem} belong to different Euclidean spaces, the equivalence between these minimization problems should be clarified. For our convenience, put $I=\{1,\dots,m\}$ and $J=\{1,\dots,k\}$. 

\medskip
Given a vector $\bar {x}=(\bar {x}^1,\dots,\bar {x}^k)\in\mathbb R^{n\times k}$, we inductively construct~$k$ subsets $A^1,\dots,A^k$ of $A$ in the following way. Put $A^0=\emptyset$ and
\begin{equation}\label{natural_clustering}
	A^j=\left\{{a}^i\in A\setminus \left(\bigcup_{p=0}^{j-1} A^p\right)\; \mid\; \|{a}^i-\bar {x}^j\|=\displaystyle\min_{q\in J} \|{a}^i-\bar {x}^q\|\right\}
\end{equation}
for $j\in J$. This means that, for every $i\in I$, the data point ${a}^i$ belongs to the cluster $A^j$ if and only if the distance $\|{a}^i-\bar {x}^j\|$ is the minimal one among the distances $\|{a}^i-\bar {x}^q\|$, $q\in J$, and ${a}^i$ does not belong to any  cluster $A^p$ with $1\leq p\leq j-1$. We will call this family $\{A^1,\dots,A^k\}$ \textit{the natural clustering associated with $\bar {x}$}.  
  \begin{Definition}\label{def_attraction_set}
  	{\rm  Let $\bar{x}=(\bar {x}^1,\dots,\bar {x}^k)\in\mathbb R^{n\times k}$. We say that the component $\bar {x}^j$ of $\bar{x}$ is \textit{attractive} with respect to the data set $A$ if the set
  		$$A[\bar {x}^j]:=\left\{{a}^i\in A\; \mid\; \|{a}^i-\bar {x}^j\|=\displaystyle\min_{q\in J} \|{a}^i-\bar {x}^q\|\right\}$$ is nonempty. The latter is called the \textit{attraction set} of $\bar {x}^j$.}
  \end{Definition}
  
  Clearly, the cluster $A^j$ in \eqref{natural_clustering} can be represented as follows:
  $$A^j=A[\bar x^j]\setminus \left(\bigcup_{p=1}^{j-1} A^p\right).$$
  
\begin{Proposition}\label{equivalence_prop} {\rm (\cite[Proposition~1]{CYY_p1})}If $(\bar{x},\bar{\alpha})$ is a solution of \eqref{basic_clustering_problem}, then $\bar{x}$ is a solution of \eqref{DC_clustering_problem}. Conversely, if  $\bar{x}$ is a solution of \eqref{DC_clustering_problem}, then the natural clustering defined by \eqref{natural_clustering} yields an incident matrix $\bar{\alpha}$ such that  $(\bar{x},\bar{\alpha})$ is a solution of \eqref{basic_clustering_problem}.
\end{Proposition}

\begin{Remark}\label{pairwise distinct data points} {\rm In practical measures, some data points can coincide. Naturally, if $a^{i_1}=a^{i_2}$, $i_1\neq i_2$, then $a^{i_1}$ and $a^{i_2}$ must belong to the same cluster. Procedure \eqref{natural_clustering} guarantees the fulfillment of this natural requirement. By grouping identical data points and choosing from each group a unique representative, we obtain a new data set having pairwise distinct data points. Thus, \textit{there is no loss of generality in assuming that ${a}^1,...,{a}^m$ are pairwise distinct points}.}
\end{Remark}
	
\begin{Theorem}\label{thm_basic_facts} {\rm (\cite[Theorem~1]{CYY_p1})}
		Both problems \eqref{basic_clustering_problem}, \eqref{DC_clustering_problem} have solutions.  If ${a}^1,...,{a}^m$ are pairwise distinct points, then the solution sets are finite. Moreover, in that case, if $\bar{x}=(\bar {x}^1,\dots,\bar {x}^k)\in\mathbb R^{n\times k}$ is  a global solution of \eqref{DC_clustering_problem}, then the attraction set $A[\bar {x}^j]$ is nonempty for every $j\in J$ and one has
	\begin{equation}\label{solution_components}
		\bar {x}^j=\frac{1}{|I(j)|}\displaystyle\sum_{i\in I(j)} {a}^i,
	\end{equation} where $I(j):=\left\{i\in I\mid {a}^i\in A[\bar {x}^j]\right\}$ with $|\Omega|$ denoting the number of elements of $\Omega$.
\end{Theorem}
 
 \begin{Proposition}\label{thm_basic_facts(1)} {\rm (\cite[Proposition~3]{CYY_p1})} If $\bar{x}=(\bar x^1,\dots,\bar x^k)\in\mathbb R^{n\times k}$ is a global solution of \eqref{DC_clustering_problem}, then the components of $\bar x$ are pairwise distinct, i.e., $\bar x^{j_1}\neq\bar x^{j_2}$ whenever $j_2\neq j_1$.
 \end{Proposition}
 
Formula \eqref{solution_components} is effective for computing certain components of any given \textit{local solution} of  \eqref{DC_clustering_problem}. The precise statement of this result is as follows.

\begin{Theorem}\label{thm_local_solutions} {\rm (\cite[Theorem~2]{CYY_p1})}
	If $\bar{x}=(\bar {x}^1,\dots,\bar {x}^k)\in\mathbb R^{n\times k}$ is a local solution of \eqref{DC_clustering_problem}, then \eqref{solution_components} is valid for all $j\in J$ whose index set $I(j)$ is nonempty, i.e., the component $\bar {x}^j$ of $\bar{x}$ is attractive w.r.t. the data set $A$. 
	
\end{Theorem}
 
 Given an  arbitrary nonempty subset $\Omega=\{a^{i_1},\dots,a^{i_r}\}\subset A$, we denote the barycenter of $\Omega$ by $b_\Omega$, that is,  $b_\Omega=\displaystyle\frac{1}{r}\sum_{l=1}^{r}a^{i_l}$. The set of all the points $b_\Omega$, where $\Omega$ is a nonempty subset of $A$, is denoted by ${\mathcal B}$.
 
 \medskip
 To proceed furthermore, we need to introduce the following condition on the local solution $x$.

\medskip
{\bf (C1)} \textit{The components of $x$ are pairwise distinct, i.e., $x^{j_1}\neq x^{j_2}$ whenever $j_2\neq j_1$.}

  \begin{Definition}\label{nontrivial_local_solution}
  	{\rm  A local solution $x = (x^1,..., x^k)$ of  \eqref{DC_clustering_problem} that satisfies {\bf (C1)} is called a \textit{nontrivial local solution}.}
  \end{Definition}
  
  As observed in \cite[Remark~4]{CYY_p1}, Proposition~\ref{thm_basic_facts(1)} shows that \textit{every global solution of \eqref{DC_clustering_problem} is a nontrivial local solution}.

\medskip
The following fundamental facts have the origin in \cite[pp.~346]{Ordin_Bagirov_2015}. 

\begin{Theorem}\label{necessary_optimality_condition} {\rm (\cite[Theorem~3]{CYY_p1})} If $x = (x^1,..., x^k)$ is a nontrivial local solution of \eqref{DC_clustering_problem} then, for any $i\in I$, $|J_i(x)|=1$. Moreover, for every $j\in J$ such that the attraction set $A[x^j]$ of $x^j$ is nonempty, one has
	\begin{equation}\label{local_solution_components}
	x^j=\frac{1}{|I(j)|}\displaystyle\sum_{i\in I(j)} a^i,
	\end{equation} where $I(j)=\left\{i\in I\mid a^i\in A[x^j]\right\}$. For any $j\in J$ with $A[x^j]=\emptyset$, one has 
	\begin{equation}\label{abnormal_components}
	x^j\notin {\mathcal A}[x],
	\end{equation} where $ {\mathcal A}[x]$ is the union of the balls $\bar B(a^p,\|a^p-x^q\|)$ with $p\in I$, $q\in J$ satisfying $p\in I(q)$.
\end{Theorem}
 
Roughly speaking, the necessary optimality condition given in the above theorem is a sufficient one. Therefore, in combination with Theorem \ref{necessary_optimality_condition}, the next statement gives a complete description of the nontrivial local solutions of \eqref{DC_clustering_problem}.

 \begin{Theorem}\label{sufficient_optimality_condition} {\rm (\cite[Theorem~4]{CYY_p1})}  Suppose that $x = (x^1,..., x^k) \in\mathbb R^{n\times k}$ satisfies condition  {\bf (C1)}. If \eqref{local_solution_components} is valid for every $j\in J$ such that the attraction set $A[x^j]\neq\emptyset$ and if \eqref{abnormal_components} is fulfilled for any $j\in J$ with $A[x^j]=\emptyset$,  then $x$ is a nontrivial local solution of \eqref{DC_clustering_problem}.
 \end{Theorem}

We complete this section by recalling the \textit{$k$-means clustering algorithm}, which is denoted by KM. This presentation of KM was given in \cite{CYY_p1}. For a detailed illustrative example, the interested reader is referred to \cite[Example~1]{CYY_p1}. 

\medskip
\textbf{The $k$-means Algorithm (KM)}: Despite its ineffectiveness, the $k$-means clustering algorithm (see, e.g., \cite[pp.~89--90]{Aggarwal 2014}, \cite{Jain 2010}, \cite[pp.~263--266]{Kantardzic 2011}, and \cite{MacQueen 1967}) is one of the most popular solution methods for \eqref{DC_clustering_problem}. One starts with selecting $k$ points $ {x}^1,\dots,{x}^k$ in $\mathbb R^n$ as the initial centroids. Then one inductively constructs~$k$ subsets $A^1,\dots,A^k$ of the data set $A$ by putting $A^0=\emptyset$ and using the rule \eqref{natural_clustering}, where ${x}^j$ plays the role of $\bar {x}^j$ for all $j\in J$. This means that $\{A^1,\dots,A^k\}$ is the natural clustering associated with $x=({x}^1,\dots,{x}^k)$. Once the clusters are formed, for each $j\in J$, if $A^j\neq\emptyset$ then the centroid ${x}^j$ is updated by the rule 
\begin{equation}\label{new_centroids}
{x}^j\ \; \leftarrow\ \; \widetilde {x}^j:=\frac{1}{|I(A^j)|}\, \displaystyle\sum_{i\in I(A^j)} {a}^i
\end{equation} with $I(A^j):=\left\{i\in I\mid {a}^i\in A^j\right\}$; and ${x}^j$ does not change otherwise. The algorithm iteratively repeats the procedure until the centroid system $\{{x}^1,\dots,{x}^k\}$ is stable, i.e., $\widetilde {x}^j={x}^j$ for all $j\in J$ with  $A^j\neq\emptyset$. The computation procedure is described as follows. 

\medskip
\noindent
\rule[0.05cm]{15.05cm}{0.01cm}\\
\textbf{Input}: The data set $A = \{{a}^1,...,{a}^m\}$ and a constant $\varepsilon\geq 0$ (tolerance).\\
\textbf{Output}: The set of $k$ centroids $\{{x}^1,...,{x}^k\}$.\\
\textit{Step 1}. Select initial centroids ${x}^j \in \mathbb{R}^n$ for all $j\in J$.\\
\textit{Step 2}. Compute $ \alpha_i=\min\{\|{a}^i-{x}^j\| \mid j\in J\}$ for all $i\in I$.\\
\textit{Step 3}.  Form the clusters $A^1,\dots,A^k$:\\
- Find the attraction sets 
$$A[{x}^j]=\left\{{a}^i\in A\; \mid\; \|{a}^i-{x}^j\|= \alpha_i\right\}\quad\; (j\in J);$$ 
- Set $A^1=A[{x}^1]$ and 
\begin{equation*}\label{natural_clusters}
A^j=A[{x}^j]\setminus \left(\bigcup_{p=1}^{j-1} A^p\right)\quad\; (j=2,\dots,k).
\end{equation*}
\textit{Step 4}. Update the centroids ${x}^j$ satisfying $A^j\neq\emptyset$ by the rule \eqref{new_centroids}, keeping other centroids unchanged.\\
\textit{Step 5}. Check the convergence condition: 
If $\|\widetilde{x}^j-{x}^j\|\leq\varepsilon$ for all $j\in J$ with $A^j\neq\emptyset$ then stop, else go to \textit{Step 2}.\\
\rule[0.05cm]{15.05cm}{0.01cm}

\section{Ordin-Bagirov's Clustering Algorithm}\label{incremental clustering algorithms}
\setcounter{equation}{0} 

In this section, we will describe the incremental heuristic algorithm of Ordin and Bagirov \cite[pp.~349--353]{Ordin_Bagirov_2015} and establish some properties of the algorithm.

\subsection{Basic constructions}\label{Basic constructions}

Let $\ell$ be an index with $1\leq\ell\leq k-1$ and let $\bar x=(\bar x^1,...,\bar x^\ell)$ be an approximate solution of problem~\eqref{DC_clustering_problem} where $k$ is replaced by $\ell$. So, $\bar x=(\bar x^1,...,\bar x^\ell)$ solves approximately the problem
\begin{equation}\label{auxDC_clustering_problem}
\min\Big\{f_\ell(x):=\frac{1}{m}\sum_{i=1}^{m}\left(\min_{j=1,\dots,\ell} \|{a}^i-{x}^j\|^2\right)\mid x=({x}^1,\dots,{x}^\ell)\in\mathbb R^{n\times\ell}\Big\}.
\end{equation} 
Applying the natural clustering procedure described in \eqref{natural_clustering} to the centroid system $\big\{\bar x^1,...,\bar x^\ell\big\}$, one divides $A$ into $\ell$ clusters with the centers $\bar x^1,...,\bar x^\ell$.  For every $i\in I$, put 
\begin{eqnarray}\label{distance_function}
d_\ell(a^i) = \min\left\{\|\bar x^1-a^i\|^2,...,\|\bar x^\ell-a^i\|^2\right\}.
\end{eqnarray} 
The formula $g(y) = f_{\ell+1}(\bar x^1,...,\bar x^\ell,y)$ where, in accordance with \eqref{auxDC_clustering_problem}, $$f_{\ell+1}(x)=\frac{1}{m}\sum_{i=1}^{m}\left(\min_{j=1,\dots,\ell+1} \|{a}^i-{x}^j\|^2\right)\quad \forall  x=({x}^1,\dots,{x}^\ell,{x}^{\ell+1})\in\mathbb R^{n\times (\ell+1)},$$ defines our \textit{auxiliary cluster function} $g:\mathbb R^n\to \mathbb R$.  From \eqref{distance_function} it follows that
\begin{eqnarray}\label{g(y)}
g(y)= \frac{1}{m}\sum_{i=1}^{m}\min\left\{d_\ell(a^i),\, \|y-a^i\|^2\right\}.\end{eqnarray} 
The problem
\begin{eqnarray}\label{aux clus prob}
\min\big\{g(y)\,\mid\, y \in \mathbb R^n\big\}
\end{eqnarray} 
is called \textit{the auxiliary clustering problem}. For each $i\in I$, one has
$$\min\left\{d_\ell(a^i),\, \|y-a^i\|^2\right\}=\big[d_\ell(a^i) +\|y-a^i\|^2\big]-\max\big\{d_\ell(a^i),\|y-a^i\|^2\big\}.$$
So, the objective function of~\eqref{aux clus prob} can be represented as $g(y) = g^1(y)-g^2(y),$ where 
\begin{eqnarray}\label{dc function} g^1(y) =  \frac{1}{m}\sum_{i=1}^{m}d_\ell(a^i) + \frac{1}{m}\sum_{i=1}^{m}\|y-a^i\|^2\end{eqnarray}  is a smooth convex function and
\begin{eqnarray}\label{dc function2} g^2(y) =  \frac{1}{m}\sum_{i=1}^{m}\max\big\{d_\ell(a^i),\|y-a^i\|^2\big\}.
\end{eqnarray} is a nonsmooth convex function. Consider the open set
\begin{equation}\label{Y_1} Y_1 := \displaystyle\bigcup_{i\in I} B\big(a^i,d_{\ell}(a^i)\big)\\
\, =\big\{y\in\mathbb R^n\,\mid\,\exists i\in I\text{ with } \|y-a^i\|^2 < d_{\ell}(a^i)\big\},\end{equation}  which is the finite union of certain open balls with the centers $a^i$ ($i\in I$), and put
$$Y_2 := \mathbb R^n\setminus Y_1=\big\{y\in\mathbb R^n\,\mid\,\|y-a^i\|^2 \geq d_{\ell}(a^i),\ \forall i\in I\big\}.$$
One sees that all the points $\bar x^1,...,\bar x^\ell$ are contained in $Y_2$. Since $\ell<k\leq m$ and the data points $a^1,\dots,a^m$ are pairwise distinct, there must exist at least one $i\in I$ with $d_{\ell}(a^i)>0$ (otherwise, every data point coincides with a point from the set $\big\{\bar x^1,...,\bar x^\ell\big\}$, which is impossible). Hence $Y_1\neq\emptyset$. By \eqref{dc function} and \eqref{dc function2}, we have \begin{equation*}\label{strict_ineq} g(y) < \frac{1}{m}\sum_{i=1}^{m}d_{\ell}(a^i)\quad\forall y\in Y_1\end{equation*} 
and
\begin{equation*}\label{equality}g(y) = \frac{1}{m}\sum_{i=1}^{m}d_{\ell}(a^i)\quad\forall y\in Y_2.\end{equation*}  Therefore, any \textit{iteration process} for solving \eqref{aux clus prob} should start with a point $y^0\in Y_1$. 

To find an approximate solution of~\eqref{DC_clustering_problem} where $k$ is replaced by $\ell+1$, i.e., the problem
\begin{equation}\label{auxDC_clustering_problem2}
\min\Big\{f_{\ell+1}(x):=\frac{1}{m}\sum_{i=1}^{m}\left(\min_{j=1,\dots,\ell+1} \|{a}^i-{x}^j\|^2\right)\mid x=({x}^1,\dots,{x}^{\ell+1})\in\mathbb R^{n\times (\ell+1)}\Big\},
\end{equation} 
 we can use the following procedure \cite[pp.~349--351]{Ordin_Bagirov_2015}. Fixing any $y\in Y_1$, one divides the data set $A$ into two disjoint subsets
\begin{equation*}\label{A1(y)}
A_1(y) := \{a^i \in A\mid\|y-a^i\|^2 < d_{\ell}(a^i)\}
\end{equation*} 
and
\begin{equation*}\label{A2(y)}
A_2(y) := \{a^i \in A\mid\|y-a^i\|^2 \geq d_{\ell}(a^i)\}.
\end{equation*} 
Clearly, $A_1(y)$ consists of all the data points standing closer to $y$ than to their cluster centers.  Since $y \in Y_1$, the set $A_1(y)$ is nonempty. Note that \begin{eqnarray}\label{g(y)_via_partition} g(y) = \frac{1}{m} \bigg(\sum_{a^i\in A_1(y)}\|y-a^i\|^2+\sum_{a^i\in A_2(y)}d_{\ell}(a^i)\bigg).\end{eqnarray}
Put $z_{\ell+1}(y)=f_\ell(\bar x)-g(y)$. Since $f_\ell(\bar x)=f_\ell(\bar x^1,...,\bar x^\ell)$  and  $g(y)= f_{\ell+1}(\bar x^1,...,\bar x^\ell,y)$, the quantity $z_{\ell+1}(y)>0$ expresses \textit{the decrease of the minimum sum-of-squares clustering criterion} when one replaces the current centroid system $\big\{\bar x^1,...,\bar x^\ell\big\}$ with~$\ell$ centers by the new one $\big\{\bar x^1,...,\bar x^\ell,y\big\}$ with~$\ell+1$ centers.  Thanks to the formula $$f_\ell(\bar x)= \frac{1}{m} \sum_{a^i\in A}d_{\ell}(a^i)$$ and \eqref{g(y)_via_partition}, one has the representation
$$z_{\ell+1}(y) = \frac{1}{m}\sum_{a^i\in A_1(y)}\big(d_{\ell}(a^i)-\|y-a^i\|^2\big),$$
which can be rewritten as
\begin{eqnarray}\label{zly}
z_{\ell+1}(y) = \frac{1}{m}\sum_{i\in I}\max\big\{0,d_{\ell}(a^i)-\|y-a^i\|^2\big\}.
\end{eqnarray}

Further operations depend greatly on the data points belonging to $Y_1$. It is easy to show that $a \in A\cap Y_1$ if and only if $a\in A$ and $a\notin\{\bar x^1,...,\bar x^\ell\}$. For every point $y = a \in A\cap Y_1$, one computes $z_{\ell+1}(a)$ by \eqref{zly}. Then, one finds the value
\begin{eqnarray}\label{z1max}
z_{\rm max}^1 := \max\big\{z_{\ell+1}(a) \;:\; a\in A\cap Y_1\big\}.
\end{eqnarray} 
 
 The selection of `good' starting points to solve \eqref{auxDC_clustering_problem2} is controlled by two parameters: $\gamma_1\in [0,1]$ and $\gamma_2\in [0,1]$. The role of each of them will be explained later on. Since the choice of these parameters can be made from the computational experience of applying the algorithm in question, the authors of \cite{Ordin_Bagirov_2015} call their algorithm \textit{heuristic}.  
 
 Using $\gamma_1$, one can find the set
\begin{eqnarray}\label{A1bar}
\bar A_1:= \{ a\in A\cap Y_1\mid z_{\ell+1}(a)\geq\gamma_1 z_{\rm max}^1\}.
\end{eqnarray} 
For $\gamma_1 = 0$, one has $\bar A_1 = A\cap Y_1$, i.e., $\bar A_1$ consists of \textit{all} the data points belonging to $Y_1$.  In contrast, for $\gamma_1 = 1$, the set $\bar A_1$ just consists of the data points yielding the \textit{largest decrease} $z_{\rm max}^1$. (As noted by Ordin and Bagirov \cite{Ordin_Bagirov_2015}, the \textit{global $k$-means algorithm} in \cite{LVV_2003} uses one of such data points for finding a $({\ell+1})$-th centroid.) \textit{Thus, $\gamma_1 $ represents the tolerance in choosing appropriate points from $A\cap Y_1$.} For each $a\in \bar A_1$, one finds the set $A_1(a)$ and computes its barycenter, which is denoted by $c(a)$. Then, one replaces $a$ by $c(a)$,  because  $c(a)$ represents the set  $A_1(a)$ better than $a$. Since $g(c(a))\leq g(a)<f_\ell(\bar x),$ one must have $c(a)\in Y_1$. Put 
\begin{eqnarray}\label{A2bar} \bar A_2=\{c(a)\; :\; a\in \bar A_1\}.\end{eqnarray} 
For each $c\in \bar A_2$, one computes the value $z_{\ell+1}(c)$ by using \eqref{zly}. Then, we find
\begin{eqnarray}\label{z2max}
z_{\rm max}^2:= \max\big\{z_{\ell+1}(c)\;:\; c\in\bar A_2\big\}.
\end{eqnarray} 
Clearly, $z_{\rm max}^2$ is the largest decrease among the values $f_{\ell+1}(\bar x^1,...,\bar x^{\ell},c)$, where $c\in \bar A_2$, in comparison with the value $f_{\ell}(\bar x)$. 

 Using $\gamma_2$, one computes 
\begin{eqnarray}\label{A3bar}
\bar A_3 = \big\{ c \in \bar A_2\mid z_{\ell+1}(c)\geq\gamma_2 z_{\rm max}^2\big\},
\end{eqnarray} which is the 
For $\gamma_2 = 0$, one has $\bar A_3 = \bar A_2$. For $\gamma_2 = 1$,  one sees that $\bar A_3$ just contains the barycenters $c\in \bar A_2$ with the largest decrease of the objective function $g(y) = f_{\ell+1}(\bar x^1,...,\bar x^\ell,y)$ of~\eqref{aux clus prob}. (As noted in \cite[p.~315]{Ordin_Bagirov_2015}, for $\gamma_1 = 0$ and $\gamma_2 = 1$, one recovers the selection of a `good' starting point in the\textit{ modified global $k$-means algorithm} suggested by Bargirov in \cite{Bagirov_2008}.) \textit{Thus, $\gamma_2$ represents the tolerance in selecting appropriate points from $\bar A_2$.}  The set 
\begin{eqnarray}\label{Omega}\Omega:=\big\{(\bar x^1,...,\bar x^\ell,c) \mid c\in\bar A_3 \big\}\end{eqnarray} contains the `good' starting points to solve \eqref{auxDC_clustering_problem2}.

\subsection{Version 1 of Ordin-Bagirov's algorithm}

On the basis of the set~$\Omega$ in \eqref{Omega}, the computation of a set of starting points to solve problem~\eqref{auxDC_clustering_problem2} is controlled by a parameter $\gamma_3\in [1,\infty)$. One applies the $k$-means algorithm to problem \eqref{auxDC_clustering_problem2} for each initial centroid system $(\bar x^1,...,\bar x^\ell,c)\in \Omega$. In result, one obtains a set of vectors $x=(x^1,\dots,x^{\ell+1})$ from $\mathbb R^{n\times (\ell+1)}$. Denote by $\bar A_4$ the set of the components $x^{\ell+1}$ of these vectors. Then, one computes the number 
\begin{eqnarray}\label{fbarlmin_OBA}
f_{\ell+1}^{\rm min} := \min\big\{g(y)\mid y\in\bar A_4\big\}.
\end{eqnarray} 
Using $\gamma_3$, one finds the set
\begin{eqnarray}\label{A5bar_OBA}
\bar A_5 = \big\{y\in\bar A_4\mid g(y)\leq \gamma_3 f_{\ell+1}^{\rm min}\big\}.
\end{eqnarray} 
For $\gamma_3 = 1$, one sees that $\bar A_5$ contains all the points $x\in\bar A_4$ at which the function $f_{\ell+1}(x)$ attains its minimum value. In contrast, if $\gamma_3$ is large enough, then $\bar A_5 = \bar A_4$. \textit{Thus, $\gamma_3$ represents the tolerance in choosing appropriate points from $\bar A_4$.}  To solve problem~\eqref{auxDC_clustering_problem2}, one will use the points from $\bar A_5$.

\medskip
The process of finding starting points is summarized as follows.

\medskip
\noindent 
\textbf{Procedure 1} \textbf{(for finding starting points)}\\
\rule[0.05cm]{15.2cm}{0.01cm}\\
\textbf{Input:} An approximate solution $\bar x=(\bar x^1,...,\bar  x^{\ell})$ of problem \eqref{auxDC_clustering_problem}, $\ell\geq 1$.\\
\textbf{Output:} A set  $\bar A_5$ of starting points to solve problem~\eqref{auxDC_clustering_problem2}.\\
\textit{Step 1}. Select three control parameters: $\gamma_1\in [0,1],\ \gamma_2\in [0,1],\ \gamma_3 \in [1,\infty)$.\\
\textit{Step 2}. Compute $z_{\rm max}^1$ by \eqref{z1max} and the set $\bar A_1$ by \eqref{A1bar}.\\
\textit{Step 3}. Compute the set $\bar A_2$ by \eqref{A2bar}, $z_{\rm max}^2$ by \eqref{z2max}, and the set $\bar A_3$ by \eqref{A3bar}.\\
\textit{Step 4}. Using \eqref{Omega}, form the set $\Omega$.\\
\textit{Step 5}. Apply the $k$-means algorithm to problem~\eqref{auxDC_clustering_problem2} for each initial centroid system $(\bar x^1,...,\bar x^\ell,c)\in \Omega$ to get the set $\bar A_4$.\\
\textit{Step 6}. Compute the value $f_{\ell+1}^{\rm min}$ by \eqref{fbarlmin_OBA}. \\ 
\textit{Step 7}. Form the set $\bar A_5$ by \eqref{A5bar_OBA}.\\
\rule[0.05cm]{15.2cm}{0.01cm}

\medskip
Now we are able to present the original version of Ordin-Bagirov's algorithm \cite[Algorithm~2, p.~352]{Ordin_Bagirov_2015} for solving problem \eqref{DC_clustering_problem}.

\medskip
\noindent\textbf{Algorithm 1} \textbf{(Ordin-Bagirov's Algorithm, Version 1)}\\
\rule[0.05cm]{15.2cm}{0.01cm}\\
\textbf{Input}: The data set $A = \{a^1,\dots,a^m\}.$\\
\textbf{Output}: A centroid system $\{\bar x^1,\dots,\bar x^k\}$.\\
\textit{Step 1}. Compute the barycenter $a^0=\displaystyle\frac{1}{m}\sum_{i=1}^{m}a^{i}$ of the data set $A$, put $\bar x^1=a^0$, and set $\ell=1$.\\
\textit{Step 2}. If $\ell=k$, then stop. Problem \eqref{DC_clustering_problem} has been solved.\\
\textit{Step 3}. Apply \textbf{Procedure 1} to the set $\bar A_5$ of starting points.\\
\textit{Step 4}.  For each $\bar y\in\bar A_5$, apply the $k$-means algorithm to problem~\eqref{auxDC_clustering_problem2} with the starting point $(\bar x^1, . . . , \bar x^\ell, \bar y)$ to find an approximate solution $x=(x^1,\dots,x^{\ell+1})$. Denote by $\bar A_6$ the set of these solutions.\\
\textit{Step 5}. Select a point $\hat x=(\hat x^1,\dots,\hat x^{\ell+1})$ from  $\bar A_6$ satisfying
\begin{equation}\label{min_value_OBA}
f_{\ell+1}(\hat x)=\min\big\{f_{\ell+1}(x)\mid x\in \bar A_6\big\}.
\end{equation}
Define $\bar x^j:=\hat x^j,\ j=1,\dots,\ell+1$. Set $\ell:=\ell+1$ and go to \textit{Step 2}.\\
\rule[0.05cm]{15.2cm}{0.01cm}

\medskip
Depending on the sizes of the data sets, the following rule to choose the control parameters triple $\gamma=(\gamma_1,\gamma_2,\gamma_3)$ can be used \cite[p.~352]{Ordin_Bagirov_2015}:

$\bullet$ For small data sets (with
the number of data points $m\leq 200$), one can choose $\gamma=(0.3,0.3,3)$; 

$\bullet$  For medium size data sets ($200 < m \leq 6000$), one can choose $\gamma=(0.5,0.8,1.5)$, or $\gamma=(0.5,0.9,1.5)$; 

$\bullet$ For large data sets (with $m > 6000$), one can choose $\gamma=(0.85,0.99,1.1)$, or $\gamma=(0.9,0.99,1.1)$.

\medskip Going back to Procedure 1 and Algorithm 1, we have the following remarks.
When one applies the $k$-means algorithm to problem \eqref{auxDC_clustering_problem2} for an initial centroid system $(\bar x^1,...,\bar x^\ell,c)\in \Omega$ to get the new centroid system $x=(x^1,\dots,x^{\ell+1})$ and put $\bar y=x^{\ell+1}$, then $\bar y$ is good just in the combination with the centroids $x^1,\dots,x^{\ell}$. If one combines $\bar y$ with the given centroids $\bar x^1,...,\bar x^\ell$, as it is done in Step 4 of the above algorithm, then it may happen that $f_{\ell+1}(x^1,\dots,x^{\ell},\bar y)<f_{\ell+1}(\bar x^1,...,\bar x^\ell,\bar y).$ If so, one wastes the available centroid system $(x^1,\dots,x^{\ell},\bar y)$ with $\bar y\in\bar A_5$. And the application of the $k$-means algorithm to problem~\eqref{auxDC_clustering_problem2} with the starting point $(\bar x^1, . . . , \bar x^\ell, \bar y)$ to find an approximate solution $x=(x^1,\dots,x^{\ell+1})$, as suggested in Step 4 of the above algorithm, is not very suitable. These remarks lead us to proposing Version 2 of  Ordin-Bagirov's algorithm, which is simpler than the original version.

\subsection{Version 2 of  Ordin-Bagirov's algorithm}

The computation of an approximate solution of problem~\eqref{auxDC_clustering_problem2}  on the basis of the set~$\Omega$ in \eqref{Omega} is controlled by a parameter $\gamma_3\in [1,\infty)$. One applies the $k$-means algorithm to problem \eqref{auxDC_clustering_problem2} for each initial centroid system $(\bar x^1,...,\bar x^\ell,c)\in \Omega$. In result, one obtains a set of points $x=(x^1,\dots,x^{\ell+1})$ from $\mathbb R^{n\times (\ell+1)}$, which is denoted by $\widetilde A_4$. Then, one computes the number 
\begin{eqnarray}\label{fbarlmin}
\widetilde f_{\ell+1}^{\rm min} := \min\big\{f_{\ell+1}(x)\mid x\in\widetilde A_4\big\}.
\end{eqnarray} Clearly, $\widetilde f_{\ell+1}^{\rm min}$ equals to the value $f_{\ell+1}^{\rm min}$ defined by \eqref{fbarlmin_OBA}.

Using $\gamma_3$, one finds the set
\begin{eqnarray}\label{A5_widetilde}
\widetilde A_5 = \big\{x\in\widetilde A_4\mid f_{\ell+1}(x) \leq \gamma_3 \widetilde f_{\ell+1}^{\rm min}\big\}.
\end{eqnarray} 
For $\gamma_3 = 1$, one sees that $\widetilde A_5$ contains all the points $x\in\widetilde A_4$ at which the function $f_{\ell+1}(x)$ attains its minimum value. In contrast, if $\gamma_3$ is large enough, then $\widetilde A_5 = \widetilde A_4$. \textit{Thus, $\gamma_3$ represents the tolerance in choosing appropriate points from $\widetilde A_4$.} Selecting an arbitrary point $\hat x=(\hat x^1,\dots,\hat x^{\ell+1})$ from $\widetilde A_5$, one has an approximate solution of problem~\eqref{auxDC_clustering_problem2}.

\medskip
The above procedure for finding a new centroid system $\hat x=(\hat x^1,\dots,\hat x^{\ell+1})$ with~$\ell+1$ centers, starting from a given centroid system $\bar x=(\bar x^1,...,\bar x^\ell)$ with $\ell$ centers, can be described as follows.

\bigskip
\noindent 
\textbf{Procedure 2} \textbf{(for finding a new centroid system)}\\
\rule[0.05cm]{15.2cm}{0.01cm}\\
\textbf{Input:} An approximate solution $\bar x=(\bar x^1,...,\bar  x^{\ell})$ of problem \eqref{auxDC_clustering_problem}, $\ell\geq 1$.\\
\textbf{Output:} An approximate solution $\hat x=(\hat x^1,\dots,\hat x^{\ell+1})$ of problem~\eqref{auxDC_clustering_problem2}.\\
\textit{Step 1}. Select three control parameters: $\gamma_1\in [0,1],\ \gamma_2\in [0,1],\ \gamma_3 \in [1,\infty)$.\\
\textit{Step 2}. Compute $z_{\rm max}^1$ by \eqref{z1max} and the set $\bar A_1$ by \eqref{A1bar}.\\
\textit{Step 3}. Compute the set $\bar A_2$ by \eqref{A2bar}, $z_{\rm max}^2$ by \eqref{z2max}, and the set $\bar A_3$ by \eqref{A3bar}.\\
\textit{Step 4}. Using \eqref{Omega}, form the set $\Omega$.\\
 \textit{Step 5}. Apply the $k$-means algorithm to problem~\eqref{auxDC_clustering_problem2} for each initial centroid system $(\bar x^1,...,\bar x^\ell,c)\in \Omega$ to get the set $\widetilde A_4$ of candidates for approximate solutions of \eqref{auxDC_clustering_problem2}.\\
\textit{Step 6}. Compute the value $\widetilde f_{\ell+1}^{\rm min}$ by~\eqref{fbarlmin} and the set $\widetilde A_5$ by~\eqref{A5_widetilde}.\\ 
\textit{Step 7}. Pick a point $\hat x=(\hat x^1,\dots,\hat x^{\ell+1})$ from $\widetilde A_5$.\\
\rule[0.05cm]{15.2cm}{0.01cm}

\medskip
Now we are able to present Version 2 of Ordin-Bagirov's algorithm \cite[Algorithm~2, p.~352]{Ordin_Bagirov_2015} for solving problem \eqref{DC_clustering_problem}.

\bigskip
\noindent\textbf{Algorithm 2} \textbf{(Ordin-Bagirov's Algorithm, Version 2)}\\
\rule[0.05cm]{15.2cm}{0.01cm}\\
\textbf{Input}: The parameters $n, m, k,$ and the data set $A = \{a^1,\dots,a^m\}.$\\
\textbf{Output}:  A centroid system $\bar x=(\bar x^1,\dots,\bar x^k)$ and the corresponding clusters $A^1,...,A^k$.\\
\textit{Step 1}. Compute the barycenter $a^0=\displaystyle\frac{1}{m}\sum_{i=1}^{m}a^{i}$ of the data set $A$, put $\bar x^1=a^0$, and set $\ell=1$.\\
\textit{Step 2}. If $\ell=k$, then go to \textit{Step 5}.\\
\textit{Step 3}. Apply \textbf{Procedure 2} to compute an approximate solution $\hat x=(\hat x^1,\dots,\hat x^{\ell+1})$ of problem~\eqref{auxDC_clustering_problem2}.\\
\textit{Step 4}. Put $\bar x^j:=\hat x^j,\ j=1,\dots,\ell+1.$ Set $\ell:=\ell+1$ and go to \textit{Step 2}.\\
\textit{Step 5}. Select an element $\bar x=(\bar x^1,\dots,\bar x^k)$ from  the set
\begin{eqnarray}\label{A6_widetilde}
\widetilde A_6 := \big\{x\in\widetilde A_5\mid f_{\ell+1}(x)=\widetilde f_{\ell+1}^{\rm min}\big\}.
\end{eqnarray}
 Using the centroid system $\bar x$, apply the natural clustering procedure to partition $A$ into $k$ clusters $A^1,...,A^k$. Print $\bar x$ and $A^1,...,A^k$. Stop.\\
\rule[0.05cm]{15.2cm}{0.01cm}

\medskip
To understand the performances of Algorithms 1 and 2, let us analyze two useful numerical examples of the MSSC problem in the form \eqref{DC_clustering_problem}. For the sake of clarity and simplicity, data sets with only few data points, each has just two features, are considered.

\begin{Example}\label{Example1} {\rm Choose $n=2$, $m=3$, $k = 2$, $A = \{a^1, a^2, a^3\}$, where $a^1 = (0,0)$, $a^2 = (1,0)$, $a^3 = (0,1)$. Let $\gamma_1 = \gamma_2 = 0.3, \gamma_3 = 3.$ The barycenter of $A$ is $a^0= (\frac{1}{3},\frac{1}{3})$.
		
The implementation of Alogrithm~1 begins with computing $\bar x^1=a^0$ and setting $\ell = 1$. Since $\ell<k$, we apply Procedure~1 to compute the set $\bar A_5$. By \eqref{distance_function}, one has $d_1(a^1) = \frac{2}{9}$, $d_1(a^2) = \frac{5}{9}$, $d_1(a^3) = \frac{5}{9}$. Using \eqref{zly}, we get $z_{\ell+1}(a^1)=\frac{2}{27}$, $z_{\ell+1}(a^2) = \frac{5}{27}$, and $z_{\ell+1}(a^3) = \frac{5}{27}$. So, by \eqref{z1max} and \eqref{A1bar}, one has $z_{\rm max}^1 =  \max\{\frac{2}{27},\,\frac{5}{27},\,\frac{5}{27} \} = \frac{5}{27}$ and $\bar A_1 = A$. Since $A_1(a^i) = \{a^i\}$ for $i\in I$, one obtains $c(a^i)=a^i$ for all $i\in I$.  Therefore, by \eqref{A2bar} and \eqref{z2max}, $\bar A_2 = A$ and $$z_{\rm max}^2 = \max\Big\{\frac{2}{27},\,\frac{5}{27},\,\frac{5}{27}\Big\} =\frac{5}{27}.$$ It follows that $\bar A_3 = \{a^1, a^2, a^3\}.$ Next, one applies the $k$-means algorithm to problem~\eqref{auxDC_clustering_problem2} with initial points from the $\Omega$ defined by \eqref{Omega} to compute $\bar A_4$. Starting from $(\bar x^1,a^1)\in\Omega$, one obtains the centroid system $\{(\frac{1}{2},\frac{1}{2}),(0,0)\}$. Starting from $(\bar x^1,a^2)$ and $(\bar x^1,a^3)$, one gets, respectively, the centroid systems $\{(\frac{1}{2},\frac{1}{2}),(0,0)\}$, $\{(0,\frac{1}{2}),(1,0)\}$, and $\{(\frac{1}{2},0),(0,1)\}$. Therefore, $\bar A_4 = \{(0,0),(1,0),(0,1)\}$. By \eqref{g(y)}, we have $g(a^1) = \frac{10}{27}$, $g(a^2) = \frac{7}{27}$, and $g(a^3) = \frac{7}{27}$. So, by \eqref{fbarlmin_OBA} one obtains   $f_{\ell+1}^{\rm min} = \frac{7}{27}$. So, from \eqref{A5bar_OBA} it follows that $\bar A_5 = \{(0,0), (1,0), (0,1)\}$. Applying again the $k$-means algorithm to problem~\eqref{auxDC_clustering_problem2} with the initial points $(\bar x^1,\bar y)$, $\bar y\in \bar A_5$, one gets $$\bar A_6 =\Big\{\Big((\frac{1}{2},\frac{1}{2}),\,(0,0)\Big),\,\Big((0,\frac{1}{2}),\,(1,0)\Big),\,\Big((\frac{1}{2},0),(0,1)\Big)\Big\}.$$
The set of the values $f_{\ell+1}(x)$, $x\in\bar A_6$, is $\Big\{\frac{1}{3},\frac{1}{6},\frac{1}{6}\Big\}$. Then, there are two centroid systems in $\bar A_6$ satisfying the condition \eqref{min_value_OBA}: \begin{eqnarray}\label{result_OBA}\hat x=\Big((0,\frac{1}{2}),\,(1,0)\Big)\quad {\rm and}\quad \hat x=\Big((\frac{1}{2},0),(0,1)\Big).\end{eqnarray}  Select any one from these centroid systems and increase $\ell$ by 1. Since $\ell=2$, i.e., $\ell=k$, the computation ends. In result, one of the two centroid systems described by~\eqref{result_OBA} is found.

The implementation of Alogrithm~2 begins with putting $\bar x^1=a^0$ and setting $\ell = 1$. Since $\ell<k$, we apply Procedure 2 to compute an approximate solution $\hat x=(\hat x^1,\dots,\hat x^{\ell+1})$ of problem~\eqref{auxDC_clustering_problem2}. The sets $\bar A_1, \bar A_2$ and $\bar A_3$ are defined as in Algorithm~1. Hence, $\bar A_3 =\bar A_2=\bar A_1=A= \{a^1, a^2, a^3\}.$ Next, we apply the $k$-means algorithm to problem~\eqref{auxDC_clustering_problem2} with initial points from the set $\Omega$ defined by \eqref{Omega} to find $\widetilde A_4$. Since $\Omega=\{(\bar x^1, a^1),(\bar x^1, a^2),(\bar x^1, a^3)\}$, one gets $$\widetilde A_4 = \Big\{\big((\frac{1}{2},\frac{1}{2}),(0,0)\big),\,\big((0,\frac{1}{2}),(1,0)\big),\,\big((\frac{1}{2},0),(0,1)\big)\Big\}.$$ By \eqref{fbarlmin}, the set of the values $\widetilde f_{\ell+1}(x)$, $x\in\widetilde A_4$, is $\Big\{\frac{1}{3},\frac{1}{6},\frac{1}{6}\Big\}$. Using \eqref{fbarlmin}, one gets $\widetilde f_{\ell+1}^{\rm min} = \frac{1}{6}$. Since $\gamma_3=3$, by \eqref{A5_widetilde} we have $\widetilde A_5=\widetilde A_4$. Pick a point $\hat x=(\hat x^1,\hat x^{2})$ from $\widetilde A_5$. Put $\bar x^j:=\hat x^j,\ j=1,2.$ Set $\ell:=\ell+1$. Since $\ell=k$, we use \eqref{A6_widetilde} to form the set $$\widetilde A_6 =\Big\{\big((0,\frac{1}{2}),(1,0)\big),\,\big((\frac{1}{2},0),(0,1)\big)\Big\}.$$ Select any  element $\bar x=(\bar x^1,\bar x^2)$ from $\widetilde A_6$ and stop. In result, we get one of the two centroid systems in \eqref{result_OBA}. 
}
\end{Example}

In the above example, results of both Algorithm 1 and  Algorithm 2 belong to the global solution set of \eqref{DC_clustering_problem}, which consists of the two centroid systems in \eqref{result_OBA}. 

\medskip
We now present a modified version of Example \ref{Example1} to show that \textit{by Algorithm 1 (resp., Algorithm 2) one may not find a global solution of problem} \eqref{DC_clustering_problem}. In other words, even for a very small data set, Algorithm 1 (resp., Algorithm 2) may yield a local, non-global solution of~\eqref{DC_clustering_problem}.

\begin{Example}\label{Example2} {\rm  Choose $n=2$, $m=4$, $k = 2$, $A = \{a^1, a^2, a^3, a^4\}$, where  $a^1 = (0,0)$, $a^2 = (1,0)$,  $a^3 = (0,1)$, $a^4 = (1,1)$. Let $\gamma_1\in [0,1],\ \gamma_2\in [0,1],\ \gamma_3 \in [1,\infty)$ be chosen arbitrarily. The barycenter of $A$ is $a^0= (\frac{1}{2},\frac{1}{2})$. 
		
		To implement Algorithm 1, we put $\bar x^1=a^0$ and set $\ell=1$. By~\eqref{distance_function}, one has $d_1(a^i)=\frac{1}{2}$ for $i\in I$. Using~\eqref{zly}, we find that $z_{\ell+1}(a^i) = \frac{1}{8}$ for $i\in I$. So, by \eqref{z1max} and \eqref{A1bar}, one gets $z_{\rm max}^1 = \frac{1}{8}$ and $\bar A_1 = A$. Since $A_1(a^i) = \{a^i\}$ for $i\in I$, one has $c(a^i) = a^i$ for $i\in I$. Therefore, by \eqref{A2bar} and \eqref{z2max}, $\bar A_2 = \{a^1,\,a^2,\,a^3,\,a^4\}$ and $z_{\rm max}^2 =\frac{1}{8}.$
		So, $\bar A_3 = \big\{a^1,\,a^2,\,a^3,\,a^4\big\}.$
		Applying the $k$-means algorithm with the starting points $(\bar x^1,c)\in\Omega$, $c\in\bar A_3$, one obtains the centroid systems  $\big((\frac{2}{3},  \frac{2}{3}), (0,0)\big)$, $\big((\frac{1}{3},  \frac{2}{3}), (1,0)\big)$,  $\big((\frac{2}{3},  \frac{1}{3}), (0,1)\big)$, and $\big((\frac{1}{3},\frac{1}{3}), (1,1)\big)$. Therefore, we have $$\bar A_4 = \{(0,0),\,(1,0),\,(0,1),\,(1,1)\}.$$
		Due to \eqref{g(y)}, one has $g((0,0)) = g((0,1)) = g((1,0))=g((1,1))=\frac{3}{8}$. So, by \eqref{distance_function} one obtains $f_{\ell+1}^{\rm min} = \frac{3}{8}$. Thus, by \eqref{A5bar_OBA},  $\bar A_5 = \bar A_4$. For each $\bar y\in \bar A_5$, we apply the $k$-means algorithm with the starting point $(\bar x^1,\bar y)$ to solve~\eqref{auxDC_clustering_problem2}. In result, we get \begin{eqnarray}\label{A6bar_example2} \bar A_6 =\Big\{\big((\frac{2}{3}, \frac{2}{3}), (0,0)\big),\, \big((\frac{1}{3}, \frac{2}{3}),(1,0)\big),\,\big((\frac{2}{3},  \frac{1}{3}),(0,1)\big),\,\big((\frac{1}{3},\frac{1}{3}),(1,1)\big)\Big\}.
		\end{eqnarray}
		Since $f_{\ell+1}(x)=\frac{1}{3}$ for every $x\in\bar A_6$, to satisfy condition  \eqref{min_value_OBA}, one can select any point $\hat x=(\hat x^1,\hat x^{2})$ from $\bar A_6$. Define $\bar x^j:=\hat x^j,\ j=1,2$. Set $\ell:=\ell+1$. Since $\ell=k$, the computation is completed. Thus, Algorithm 1 yields one of the four centroid systems in \eqref{A6bar_example2}, which is a local, non-global solution of our clustering problem (see Remark \ref{Rem1} for detailed explanations).  
		
		The implementation of Alogrithm~2 begins with putting $\bar x^1=a^0$ and setting $\ell = 1$. Since $\ell<k$, we apply Procedure 2 to compute an approximate solution $\hat x=(\hat x^1,\dots,\hat x^{\ell+1})$ of problem~\eqref{auxDC_clustering_problem2}. The sets $\bar A_1, \bar A_2$ and $\bar A_3$ are defined as in Algorithm~1. Hence, $\bar A_3 =\bar A_2=\bar A_1=A= \{a^1, a^2, a^3,a^4\}.$ Next, we apply the $k$-means algorithm to problem~\eqref{auxDC_clustering_problem2} with initial points from the set $\Omega$ defined by \eqref{Omega} to find $\widetilde A_4$. Since $\Omega=\{(\bar x^1, a^1),(\bar x^1, a^2),(\bar x^1, a^3),(\bar x^1, a^4)\}$, one gets $$\widetilde A_4 =\Big\{\big((\frac{2}{3}, \frac{2}{3}), (0,0)\big),\, \big((\frac{1}{3}, \frac{2}{3}),(1,0)\big),\,\big((\frac{2}{3},  \frac{1}{3}),(0,1)\big),\,\big((\frac{1}{3},\frac{1}{3}),(1,1)\big)\Big\}.$$ By \eqref{fbarlmin}, the set of the values $\widetilde f_{\ell+1}(x)$, $x\in\widetilde A_4$, is $\Big\{\frac{1}{3},\frac{1}{3},\frac{1}{3},\frac{1}{3}\Big\}$. Using \eqref{fbarlmin}, one gets $\widetilde f_{\ell+1}^{\rm min} = \frac{1}{3}$. By \eqref{A5_widetilde}, we obtain $\widetilde A_5 = \widetilde A_4$. Set $\bar x = \hat x$, $\hat x\in\widetilde A_4$ and $\ell=2$. Since $\ell =k$, $\widetilde A_6=\widetilde A_5$. Select any centroid system from $\widetilde A_6$, e.g., $\bar x =\big((\frac{2}{3}, \frac{2}{3}), (0,0)\big)$. Applying the natural clustering procedure, one gets the clusters $A^1=\{a^2,a^3,a^4\}$, $A^2=\{a^1\}$, then stop.
	}
\end{Example}
\begin{Remark}\label{Rem1} {\rm Concerning the analysis given in Example~\ref{Example2}, observe that every centroid system in $\bar A_6$ is a nontrivial local solution of problem \eqref{DC_clustering_problem}. This assertion can be verified by \cite[Theorem~4]{CYY_p1}. The value of the objective function at these centroid systems is $\frac{1}{3}$. Consider a possible partition $A=A^1\cup A^2$, where $A^1$ and $A^2$ are disjoint nonempty subsets of $A$, then compute the barycenter $x^j$ of $A^j$ for $j=1,2$, and put $x=(x^1,x^2)$. According to Theorem 1,  Proposition 2, Remark 4, and Theorem 3 from \cite{CYY_p1}, global solutions of \eqref{DC_clustering_problem} do exist and belong to the set of those points~$x$. Hence, by symmetry, it is easy to see that the clustering problem in question has two global solutions: $\bar x = \big((\frac{1}{2},0),(\frac{1}{2},1)\big)$ and  $\hat x = \big((0, \frac{1}{2}),(1, \frac{1}{2})\big)$. As $f(\bar x)=f(\hat x)= \frac{1}{4}$, the four centroid systems in $\bar A_6$ is a global solution of \eqref{DC_clustering_problem} are all local, non-global solutions of  \eqref{DC_clustering_problem}. Similarly, the four centroid systems in $\widetilde A_6=\widetilde A_5 = \widetilde A_4$ are all local, non-global solutions of \eqref{DC_clustering_problem}.}
\end{Remark}

\begin{Remark}\label{Rem2}  {\rm In both Algorithm 1 and Algorithm 2, one starts with $\bar x^1=a^0$, where $a^0=\displaystyle\frac{1}{m}\sum_{i=1}^{m}a^{i}$ is the barycenter of the data set $A$. As it has been shown in Remark~\ref{Rem1}, for the clustering problem in Example~\ref{Example2} and for arbitrarily chosen control parameters $\gamma_1\in [0,1],\ \gamma_2\in [0,1],\ \gamma_3 \in [1,\infty)$, Algorithm 1 (resp., Algorithm 2) yields a local, non-global solution of~\eqref{DC_clustering_problem}. Anyway, if one starts with a data point $a^i$, $i\in I$, then by Algorithm 1 (resp., Algorithm 2) one can find a global solution of~\eqref{DC_clustering_problem}.} 
\end{Remark}

To proceed furthermore, we need next lemma.

\begin{Lemma}\label{k-means_lemma} Let $x=({x}^1,\dots,{x}^k)\in\mathbb R^{r\times k}$ be a centroid system, where the centroids ${x}^1,\dots,{x}^k$ are pairwise distinct. Then, after one step of applying the $k$-means Algorithm, one gets a new centroid system $\widetilde x=(\widetilde{x}^1,\dots,\widetilde{x}^k)$ with pairwise distinct centroids, i.e., $\widetilde{x}^{j_1}\neq \widetilde{x}^{j_2}$ for any $j_1,j_2\in J$ with $j_1\neq j_2$.
\end{Lemma}
\textit{Proof.}\ Let $\{A^1,\dots,A^k\}$ be the natural clustering associated with $x=({x}^1,\dots,{x}^k)$. For each $j\in J$, if $A^j\neq\emptyset$ then the centroid ${x}^j$ is updated by the rule \eqref{new_centroids}, and ${x}^j$ does not change otherwise. This means that
\begin{equation}\label{new_centroids-1}
\widetilde {x}^j=\frac{1}{|I(A^j)|}\, \displaystyle\sum_{i\in I(A^j)} {a}^i
\end{equation} if $A^j\neq\emptyset$, where $I(A^j)=\left\{i\in I\mid {a}^i\in A^j\right\}$, and $\widetilde {x}^j=x^j$ if $A^j=\emptyset$. Now,  suppose that $j_1,j_2\in J$ are such that $j_1\neq j_2$. We may assume that $j_1<j_2$.

Let $y^0:=\frac{1}{2}(x^{j_2}+x^{j_1})$ and $L:=\{y\in\mathbb R^n \mid \langle y-y^0, x^{j_2}-x^{j_1}\rangle =0 \}$. Then, any point $y\in\mathbb R^n$ having equal distances to $x^{j_1}$ and $x^{j_2}$ lies in $L$. Denote by $P_1$ (resp., $P_2$) the open half-space with the boundary $L$ that contains $x^{j_1}$ (resp., $x^{j_2}$).

If the clusters $A^{j_1}$ and $A^{j_2}$ are both nonempty, then $\widetilde{x}^{j_1}$ and $\widetilde{x}^{j_2}$ are defined by formula \eqref{new_centroids-1}. Since $\{A^1,\dots,A^k\}$ is the natural clustering associated with the centroid system $x=({x}^1,\dots,{x}^k)$ and $j_1<j_2$, one must have $A^{j_1}\subset\bar P_1$, where $\bar P_1:= P_1\cup L$ is the closure of $P_1$, while $A^{j_2}\subset P_2$. The formulas   
\begin{equation*}
\widetilde {x}^{j_1}=\frac{1}{|I(A^{j_1})|}\, \displaystyle\sum_{i\in I(A^{j_1})} {a}^i,\quad\ \widetilde {x}^{j_2}=\frac{1}{|I(A^{j_2})|}\, \displaystyle\sum_{i\in I(A^{j_2})} a^i
\end{equation*} show that $\widetilde {x}^{j_1}$ (resp., $\widetilde {x}^{j_2}$) is a convex combination of the points from $A^{j_1}$ (resp., $A^{j_2}$). Hence, by the convexity of $\bar P_1$ (resp., $P_2$), we have $\widetilde {x}^{j_1}\in \bar P_1$ (resp.,  $\widetilde {x}^{j_2}\in P_2$). Then, the property $\widetilde{x}^{j_1}\neq \widetilde{x}^{j_2}$ follows from the fact that $\bar P_1\cap P_2=\emptyset$.

If the clusters $A^{j_1}$ and $A^{j_2}$ are both empty, then $\widetilde{x}^{j_1}=x^{j_1}$ and $\widetilde{x}^{j_2}=x^{j_2}$. Since ${x}^1,\dots,{x}^k$ are pairwise distinct, we have $\widetilde{x}^{j_1}\neq \widetilde{x}^{j_2}$. 

If $A^{j_1}\neq\emptyset$ and $A^{j_2}=\emptyset$, then $\widetilde{x}^{j_1}\in\bar P_1$ and $\widetilde{x}^{j_2}=x^{j_2}\in P_2$. Since $\bar P_1\cap P_2=\emptyset$, one must have $\widetilde{x}^{j_1}\neq \widetilde{x}^{j_2}$. The situation $A^{j_1}=\emptyset$ and $A^{j_2}\neq\emptyset$ is treated similarly. 

The proof is complete. $\hfill\Box$

\medskip
Remarkable properties of Algorithm 2 are described in forthcoming theorems, where the following assumption is used:

\medskip
{\bf (C2)} \textit{The data points ${a}^1,...,{a}^m$ in the given data set $A$ are pairwise distinct.}

\medskip
Note that, given any data set, one can apply the trick suggested in Remark~\ref{pairwise distinct data points} to obtain a data set satisfying {\bf (C2)}.

\begin{Theorem}\label{OBA_Property1} Let $\ell$ be an index with $1\leq\ell\leq k-1$ and let $\bar x=(\bar x^1,...,\bar x^\ell)$ be an approximate solution of problem~\eqref{DC_clustering_problem} where $k$ is replaced by $\ell$. If {\bf (C2)}  is fulfilled and the centroids $\bar x^1,...,\bar x^\ell$ are pairwise distinct, then the centroids $\hat x^1,\dots,\hat x^{\ell+1}$ of the approximate solution $\hat x=(\hat x^1,\dots,\hat x^{\ell+1})$ of~\eqref{auxDC_clustering_problem2}, which is obtained by Procedure 2, are also pairwise distinct. 
\end{Theorem}
\textit{Proof.}\ Since $1\leq\ell\leq k-1$, $k\leq m$, and  data points ${a}^1,...,{a}^m$ in the given data set $A$ are pairwise distinct, one can find a data point $a^{i_0}\in A$, which is not contained in the set $\{\bar x^1,...,\bar x^\ell\}$. Then the set $Y_1$ defined by \eqref{Y_1} is nonempty, because $d_{\ell}(a^{i_0})>0$ (hence the open ball  $B\big(a^{i_0},d_{\ell}(a^{i_0})\big)$ is nonempty). Moreover, $A\cap Y_1\neq\emptyset$. So, from~\eqref{z1max} and~\eqref{A1bar} it follows that $\bar A_1\neq\emptyset$. Then, one easily deduce from \eqref{A2bar}--\eqref{A3bar} that the sets $\bar A_2$ and $\bar A_3$ are nonempty.

By the construction \eqref{Y_1} of $Y_1$, one has $Y_1\cap \{\bar x^1,...,\bar x^\ell\}=\emptyset$. It follows that $\bar A_1\cap \{\bar x^1,...,\bar x^\ell\}=\emptyset$. (Actually, this property has been noted before.) Since $$z_{\ell+1}(c(a))\geq z_{\ell+1}(a)>0\quad \forall a\in \bar A_1,$$ and $z_{\ell+1}(\bar x^j)=0$ for every $j\in \{1,\dots,\ell\}$, we have $\bar A_2\cap \{\bar x^1,...,\bar x^\ell\}=\emptyset$. As $\bar A_3\subset\bar A_2$, one sees that  $\bar A_3\cap \{\bar x^1,...,\bar x^\ell\}=\emptyset$. Consequently, by~\eqref{Omega}, the centroids in any centroid system  $(\bar x^1,...,\bar x^\ell,c)\in\Omega$ are pairwise distinct. Since the approximate solution $\hat x=(\hat x^1,\dots,\hat x^{\ell+1})$ of~\eqref{auxDC_clustering_problem2} is obtained from one centroid system  $(\bar x^1,...,\bar x^\ell,c)\in\Omega$ after applying finitely many steps of KM, thanks to Lemma~\ref{k-means_lemma} we can assert that  the centroids $\hat x^1,\dots,\hat x^{\ell+1}$  are pairwise distinct.
$\hfill\Box$

\begin{Theorem}\label{OBA_Property2} If {\bf (C2)} is fulfilled, then the centroids $\bar x^1,\dots,\bar x^k$ of the centroid system $\bar x=(\bar x^1,\dots,\bar x^k)$, which is obtained by Algorithm 2, are pairwise distinct.
\end{Theorem}
\textit{Proof.} To implement Alogrithm~2, we first compute the barycenter $a^0=\displaystyle\frac{1}{m}\sum_{i=1}^{m}a^{i}$ of the data set $A$, put $\bar x^1=a^0$, and set $\ell=1$. Next, we apply Procedure~2 to find an approximate solution $\hat x=(\hat x^1,\dots,\hat x^{\ell+1})$ of~\eqref{auxDC_clustering_problem2}. By Theorem~\ref{OBA_Property1}, the centroids $\hat x^1,...,\hat x^{\ell+1}$ are pairwise distinct. Since Procedure 2 is ended at Step 7 by picking any point $\hat x=(\hat x^1,\dots,\hat x^{\ell+1})$ from the set $\widetilde A_5$, which is defined by~\eqref{A5_widetilde}, Theorem~\ref{OBA_Property1} assures that every centroid system $\hat x=(\hat x^1,\dots,\hat x^{\ell+1})$ in $\widetilde A_5$ consists of pairwise distinct centroids. 

In Step 4 of  Algorithm 2, after putting  $\bar x^j=\hat x^j$ for $j=1,\dots,\ell+1$, one sets $\ell:=\ell+1$ and goes to Step 2. If $\ell<k$, then the computation continues, and one gets a  approximate solution $\hat x=(\hat x^1,\dots,\hat x^{\ell+1})$ of~\eqref{auxDC_clustering_problem2} with $\hat x^1,\dots,\hat x^{\ell+1}$ being pairwise distinct by Theorem~\ref{OBA_Property1}. If $\ell=k$, then the computation terminates by selecting an element $\bar x=(\bar x^1,\dots,\bar x^k)$ from  the set $\widetilde A_6$, which is defined by \eqref{A6_widetilde}. Since  $\widetilde A_6\subset\widetilde A_5$ and we have shown that  every centroid system in $\widetilde A_5$ consists of pairwise distinct centroids, the obtained centroids $\bar x^1,\dots,\bar x^k$ are pairwise distinct.
$\hfill\Box$

\medskip
On one hand, if $\bar x=(\bar x^1,\dots,\bar x^k)$ is global solution of~\eqref{DC_clustering_problem}, then by Proposition~\ref{thm_basic_facts(1)} we know that the centroids $\bar x^1,\dots,\bar x^k$ are pairwise distinct.  More general, by Definition~\ref{nontrivial_local_solution}, the components of any nontrivial local solution are pairwise distinct. On the other hand, according to Theorem~\ref{OBA_Property2}, Algorithm 2 yields a centroid system having pairwise distinct components. Thus, \textit{the  centroid system resulted from Algorithm 2 is a very good candidate for being a nontrivial local solution of~\eqref{DC_clustering_problem}}. As the global solutions are among the nontrivial local solutions, \textit{Theorem~\ref{OBA_Property2} reveals a nice feature of Algorithm 2}. 

\subsection{The $\varepsilon$-neighborhoods technique}\label{varepsilon_technique}
\textit{The $\varepsilon$-neighborhoods technique} \cite[pp.~869-870]{BUW_2011} (see also \cite[pp.~352--353]{Ordin_Bagirov_2015}) allows one to reduce the computation volume of Algorithm~1 (as well as that of Algorithm~2, or another incremental clustering algorithm based on the sets $\bar A_1$), when it is applied to large data sets. The procedure of removing data points from $A$ to get a smaller set $\bar A_1$ is as follows. Choose a sufficiently small number $\delta\in (0,\ell^{-1})$ (for example, $\delta=\min\big\{10^{-3},\ell^{-1}\big\}$). In the notations of Subsection \ref{Basic constructions}, let $\{A^1,\dots,A^\ell\}$ be the natural clustering associated with the centroid system  $\bar x=(\bar x^1,...,\bar x^\ell)$. For every $j\in\{1,\dots,\ell \}$, if $A^j\neq\emptyset$, then one defines
\begin{equation*}
\alpha_{j} = \frac{1}{|A^{j}|}\sum_{a\in A^{j}}\|\bar x^{j}-a\|^2
\end{equation*}
and $\beta_j=\max\big\{\|\bar x^{j}-a\|^2\, :\, a\in A^{j}\big\}$. Set $\mu_j = \dfrac{\beta_j}{\alpha_{j}}$ and observe that $\mu_j\geq 1$. Let
\begin{equation*}
A^{j}_\delta:=\big\{a\in A^{j}\,:\|\bar x^{j}-a\|^2\geq \eta_{j}\alpha_j\big\},
\end{equation*}  where  $\eta_j=1+\ell\delta(\mu_j-1)$. One has $A^{j}_\delta\neq\emptyset$. Indeed, if $\bar a\in A^j$ 
is such a data point that $\|\bar x^{j}-\bar a\|^2=\beta_j$, then $\|\bar x^{j}-a\|^2=\beta_j\geq \eta_{j}\alpha_j$; hence $\bar a\in A^{j}_\delta$. To proceed furthermore, denote by $A_\delta$ the union of all the sets $A^{j}_\delta$, where $j\in\{1,\dots,\ell\}$ is such that $A^j\neq\emptyset$. Now, instead of $\bar A_1$ given by \eqref{A1bar}, we use the set \begin{eqnarray}\label{A1bar-delta}
\bar A_{1,\delta}:= \{ a\in A_\delta\cap Y_1\mid z_{\ell+1}(a)\geq\eta_1 z_{\rm max}^1\},
\end{eqnarray} which is a subset of $\bar A_1$. In the construction $\bar A_{1,\delta}$ by \eqref{A1bar-delta}, we have removed from $A$ all the data points $a$ with $ \|\bar x^{j}-a\|^2< \eta_{j}\alpha_j$,  where $j\in\{1,\dots,\ell\}$ is such that $A^j\neq\emptyset$.


\begin{thebibliography}{99}
		\bibitem{Aggarwal 2014} Aggarwal, C.C., Reddy, C.K.: Data Clustering Algorithms and Applications. Chapman $\&$  Hall/CRC Press, Boca Raton, Florida (2014)
		
		\bibitem{Aloise 2009} Aloise, D., Deshpande, A., Hansen, P., Popat P.: NP-hardness of Euclidean sum-of-squares clustering. Mach. Learn. {\bf 75}, 245--248 (2009)
		
		\bibitem{Bagirov_2008}  Bagirov, A.M.: Modified global $k$-means algorithm for minimum sum-of-squares clustering problems. Pattern Recognit. {\bf 41}, 3192--3199 (2008)
		
		\bibitem{Bagirov_2014}  Bagirov, A.M.: An incremental DC algorithm for the minimum sum-of-squares clustering. Iranian J. Oper. Res. {\bf 5}, 1--14 (2014)
		
		\bibitem{Bagirov_Rubinov_2003} Bagirov, A.M., Rubinov, A.M., Soukhoroukova, N.V., Yearwood, J.: Unsupervised and supervised data classification via nonsmooth and global optimization. TOP {\bf 11}, 1--93 (2003)
					
		\bibitem{BUW_2011} Bagirov, A.M., Ugon, J., Webb, D.: Fast modified global $k$-means algorithm for incremental cluster construction. Pattern Recognit. {\bf 44}, 866--876 (2011) 
		
		\bibitem{Bagirov_Yearwood_2006} Bagirov, A.M., Yearwood, J.: A new nonsmooth optimization algorithm for minimum  sum-of-squares clustering problems. European J. Oper. Res. \textbf{170}, 578--596 (2006)
					
		\bibitem{Bock_1998} Bock, H.H.: Clustering and neural networks. In: Rizzi, A., Vichi, M., Bock, H.H. (eds.) Advances in Data Science and Classiﬁcation, pp. 265--277. Springer, Berlin (1998)
		
		\bibitem{Brusco 2006} Brusco, M.J.: A repetitive branch-and-bound procedure for minimum within-cluster sum of squares partitioning. Psychometrika {\bf 71}, 347--363 (2006)
		
		
		\bibitem{Costa 2017} Costa, L.R., Aloise, D., Mladenovi\'c, N.: Less is more: basic variable neighborhood search heuristic for balanced minimum sum-of-squares clustering. Inform. Sci. \textbf{415/416}, 247--253 (2017) 
		
		\bibitem{CYY_p1} Cuong, T.H., Yao, J.-C., Yen, N.D.: Qualitative properties of the minimum sum-of-squares clustering problem. Preprint (arXiv:1810.02057 [math.OC]), submitted.
					
		\bibitem{Du Merle 2000} Du Merle, O., Hansen, P., Jaumard, B., Mladenovi\'c, N.: An interior point algorithm for minimum sum of squares clustering. SIAM J. Sci. Comput. {\bf 21}, 1485--1505 (2000)
		
		\bibitem{HNCM_2005} Hansen, P., Ngai, E., Cheung, B.K., Mladenovic, N.: Analysis of global $k$-means, an incremental heuristic
		for minimum sum-of-squares clustering. J. Classif. \textbf{22}, 287--310 (2005)
		
		\bibitem{Jain 2010} Jain, A.K., Data clustering: 50 years beyond $k$-means, Pattern Recognit. Lett. {\bf 31}, 651--666 (2010)
		
		\bibitem{Kantardzic 2011} Kantardzic, M.: Data Mining Concepts, Models, Methods, and Algorithms, Second Edition. John Wiley $\&$ Sons, Hoboken, New Jersey (2011)
		
		\bibitem{Kumar 2017} Kumar, K.M., Reddy, A.R.M: An efficient $k$-means clustering filtering algorithm using density based initial cluster centers. Inform. Sci. \textbf{418/419}, 286--301 (2017)  
		
		\bibitem{Lai_Huang_2010} Lai, J.Z.C., Huang, T.-J.: Fast global $k$-means clustering using cluster membership and inequality. Pattern
		Recognit. \textbf{43}, 731--737 (2010)
		
		\bibitem{PhamDinh_LeThi_2009} Le Thi, H.A., Pham Dinh, T.: Minimum sum-of-squares clustering by DC programming and DCA. ICIC 2009, LNAI \textbf{5755}, 327--340 (2009)
			
		\bibitem{LVV_2003} Likas, A., Vlassis, N., Verbeek, J.J.: The global $k$-means clustering algorithm. Pattern Recognit. {\bf 36}, 451--461 (2003)
		
		\bibitem{MacQueen 1967} MacQueen, J.: Some methods for classification and analysis of multivariate observations. In: Proceedings of the 5th Berkeley Symposium on Mathematical Statistics and Probability, pp. 281--297 (1967)
		
		\bibitem{Mahajan_NPHard_2009} Mahajan, M., Nimbhorkar, P., Varadarajan, K.: The planar $k$-means problem is NP-hard. Theoret. Comput. Sci. \textbf{442}, 13--21 (2012)
		
		\bibitem{Ordin_Bagirov_2015}  Ordin, B., Bagirov, A.M.:  A heuristic algorithm for solving the minimum sum-of-squares clustering problems. J. Global Optim. {\bf 61}, 341--361 (2015)
		
			
		\bibitem{Peng 2005} Peng, J., Xiay, Y.: A cutting algorithm for the minimum sum-of-squared error clustering. In: Proceedings of the SIAM International Data Mining Conference (2005)
		
		\bibitem{Shereli 2005} Sherali, H.D., Desai, J.: A global optimization RLT-based approach for solving the hard clustering problem. J. Global Optim. {\bf 32}, 281--306 (2005)
				
		\bibitem{Xie_2011} Xie, J., Jiang, S., Xie, W., Gao, X.:  An efficient global $k$-means clustering algorithm. J. Comput. {\bf 6}, 271--279 (2011)	 		
	
\end{thebibliography}
\end{document}